\documentclass[11pt,twoside]{article}
\usepackage{graphics,latexsym,amssymb,amsmath}
\RequirePackage{graphics,epsfig,psfrag}
\def\abs#1{\left \vert #1 \right \vert}
\def\RR{{\bf R}} 
\def\SS{{\bf S}} 
\def\H{{\rm H}}
\def\pn{\medskip\par\noindent}
\def\Frac#1#2{{\displaystyle{{#1}\over{#2}}}}

\def\[#1\]{\begin{eqnarray}#1\end{eqnarray}}
\def\$#1\${\begin{eqnarray}#1\end{eqnarray}}

\def\phi{\varphi}

\def\eps{\varepsilon}

\def\cC{{\cal C}}

\def\pent#1#2{\pe{\frac{#1}{#2}}}

\def\pe#1{{\left \lbrack #1 \right \rbrack}} 

\def\sign#1{{\rm sign}\,\bigl( #1 \bigr)}
\newcommand{\Pf}{{\em Proof}. }

\newcommand{\EPf}{\hbox{}\hfill$\Box$\vspace{.5cm}}

\newcommand{\ZZ}{{\mbox{\bf Z}}}
\newcommand{\T}{{\mathrm{T}}}
\def\Mod#1{\,(\mathrm{mod}\,#1)}
\def\Sum{\mathop{\sum}\limits}

\def\abs#1{\left \vert #1 \right \vert}
\def\Frac#1#2{{\displaystyle{\frac{#1}{#2}}}}
\def\phi{\varphi}

\newtheorem{definition}{Definition}
\newtheorem{remark}{Remark}
\newtheorem{lemma}{Lemma}
\newtheorem{proposition}{Proposition}
\newtheorem{corollary}{Corollary}
\newtheorem{thm*}{Theorem}
\newtheorem{theorem}{Theorem}
\date{\today}
\begin{document}
\pagestyle{myheadings}
\markboth{P. -V. Koseleff, D. Pecker}{{\em Chebyshev Knots}}
\title{ Chebyshev knots }
\author{P. -V. Koseleff\footnote{Facult{\'e} de Math{\'e}matiques and Salsa-Inria project},
D. Pecker\footnote{Facult{\'e} de Math{\'e}matiques}
\medskip\\
Universit{\'e} Pierre et Marie Curie\\
4, place Jussieu, F-75252 Paris Cedex 05 \\
e-mail: {\tt\{koseleff,pecker\}@math.jussieu.fr}}
\maketitle
\begin{abstract}
A Chebyshev knot is a  knot which admits a parametrization of the
form $ x(t)=T_a(t); \  y(t)=T_b(t) ; \ z(t)= T_c(t + \phi), $ where
$a,b,c$ integers, $T_n(t)$ is the Chebyshev polynomial
of degree $n,$ and $\phi \in \RR .$ Chebyshev knots are non
compact analogues of the classical Lissajous  knots. We show that
there are infinitely many Chebyshev knots with $\phi = 0.$ We also
show that every knot is a Chebyshev knot.
\pn {\bf keywords:} {Polynomial curves, Chebyshev polynomials,  Chebyshev curves,
Lissajous knots, long knots, braids} 
\end{abstract}
\section{Introduction}
A  Lissajous knot is a knot which
admits a one-to-one parametrization of the form
$$ x=\cos (at ); \   y=\cos (bt + \phi) ; \  z=\cos (ct + \psi)
$$
where $ 0 \le t \le 2 \pi $ and where $ a, b, c$ are pairwise
coprime integers. These knots, first defined in \cite{BHJS},
have been studied by many authors: V. F. R. Jones, J. Przytycki,
C. Lamm, J. Hoste and L. Zirbel. Most  known properties of Lissajous
knots are deduced from their symmetries, which are easy to see
(see \cite{JP,La,HZ,Pr}).
\pn
On the other hand Vassiliev considered polynomial knots, i.e. non
singular polynomial  embeddings $ \RR \rightarrow \  \RR^3$
 (see \cite{Va,Sh,RS,KP1,KP2}).
\pn In this paper we study a polynomial analogue of Lissajous
knots.
It is natural to use the classical Chebyshev polynomials $T_n(t)$
instead of cosine functions to define our Chebyshev knots.
The Chebyshev  polynomials are defined by
$T_0=1, \, T_1 = t, \, T_{n+1} = 2t T_n - T_{n-1}, \, n \in \ZZ$. They satisfy
the trigonometric identity
$ \cos (n \,\theta) = T_n (\cos \,\theta)$.
\begin{definition} A knot in $\RR ^3  \subset  {\SS}^3$ is a Chebyshev
knot if it admits a one to one parametrization of the form
$$ x=T_a(t); \  y=T_b(t) ; \    z=T_c(t + \phi) $$
where $t \in \RR$, $a,b, c$ are integers
and $\phi$ is a real constant.
\end{definition}
\begin{figure}[th]
\begin{center}
{\scalebox{.4}{\includegraphics{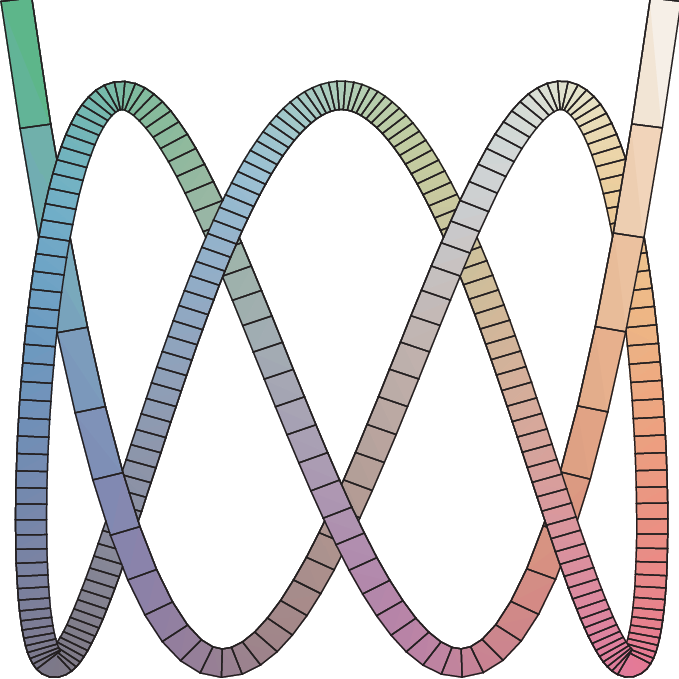}}}
\caption{The mirror image of the $7_7$ knot is Chebyshev.\label {77}}
\end{center}
\end{figure}
We begin with the study of plane Chebyshev curves which are
projections of Chebyshev knots
on the $(x,y)$-plane. We conclude this paragraph with a theorem
of Hoste and Zirbel  \cite{HZ} describing these curves in terms
of  particular braid projections. \pn Then, we study some families
of Chebyshev knots with  $ \phi = 0, $ called harmonic knots. We
prove that for $a,b$ coprime positive integers and $c= ab-a-b$, the harmonic
knot $\H(a,b,c)$ has an alternating projection  on the $(x,y)$-plane.
 We deduce that there are infinitely many harmonic
 knot types. This is similar to a theorem of C. Lamm concerning Lissajous
knots (see \cite{La}).
 We also prove that the  torus knots $\T(2,2n+1)$ are harmonic knots.
On the other hand, we observe that the symmetries of
 harmonic knots  are quite different from those of Lissajous knots.
 There are infinitely many amphicheiral harmonic knots and
 infinitely many strongly reversible harmonic knots.
The trefoil and the figure-eight knot are harmonic knots but
are not Lissajous. Some knots are both Lissajous and harmonic knots, e.g.
$5_2$ and $7_5$.
\pn
We conclude the paper with our principal result: {\em every
knot is a Chebyshev knot}. This is done by
showing first that every knot has a plane projection which is a
Chebyshev curve. Then we use
some classical results of braid theory and
a density argument based on Kronecker's theorem.
\pn
At the end we give Chebyshev diagrams of the  first 2-bridge harmonic knots.

\section{Geometry of plane Chebyshev curves}\label{cheby}
\begin{figure}[th]
\begin{center}
{\scalebox{.3}{\rotatebox{90}{\includegraphics{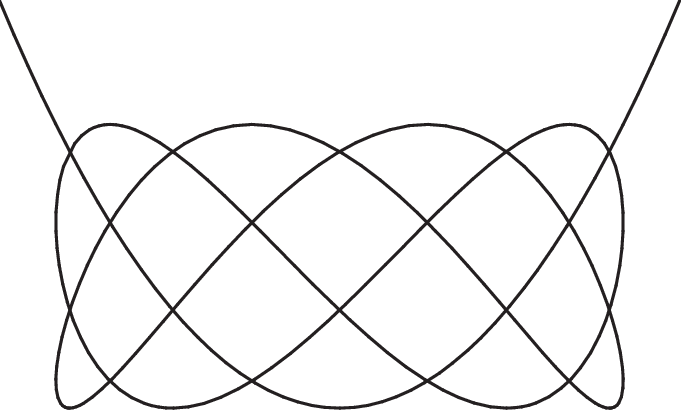}}}}\quad
{\scalebox{.3}{\rotatebox{90}{\includegraphics{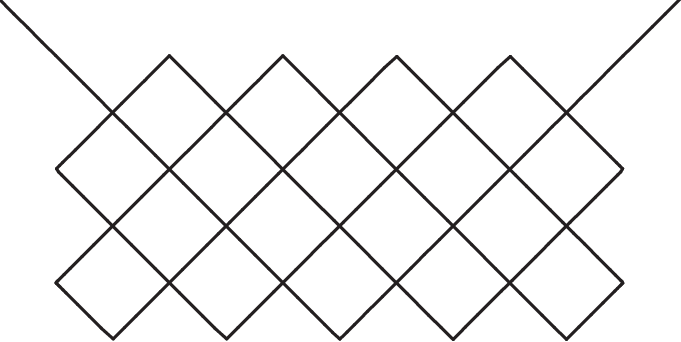}}}}
\caption{The Chebyshev curve $T_{10}(x)-T_5(y)= 0$
and its billiard picture \label{T510}}
\end{center}
\end{figure}

Chebyshev curves were defined in \cite{Fi} to replace the older
denomination of ``doubly parame\-tri\-zed Lissajous curves''.
Their double points are easier to study than those of Lissajous
curves. It will be convenient to consider also the case of
implicit Chebyshev  curves.
\begin{proposition}\label{cc}
Let $a,b$ be nonnegative integers, a being odd. The affine Chebyshev curve
${\cC}$ defined  by
$$ {\cC} \, : \quad T_b(x) - T_a(y) = 0$$
has $\frac 12{(a-1)(b-1)}$ singular points which are crossing
points. These points form two rectangular  grids contained in the
open square $Q=( -1,1)^2,$
$R=\{ (x,y) \in Q, \ T_b(x)=T_a(y)=1\},$
and
$R'= \{ (x,y) \in Q , \ T_b(x)=T_a(y)= -1 \}.$
\end{proposition}
\Pf
The singular points of ${\cC}$ are obtained for
$T'_b(x)= 0, \  T'_a(y)= 0,  \  T_b(x)=T_a(y).$
From $T_a (\cos \theta) = \cos a \theta$, we deduce that $T_a$ has
degree $a$ and $T'_a(\cos \theta) = a \Frac{\sin a \theta}{\sin
  \theta}$. $T'_a$ has $a-1$ simple roots in $(0,1)$: $y_k = \cos \Bigl(k
\Frac{\pi}{a}\Bigr)$,
$k=1,\ldots, a-1$. At these points, we have $T_a (y_k) = (-1)^k$.
$T'_b$ has $b-1$ roots in $(0,1)$: $x_1, \ldots, x_{b-1}$. For each
$x_i$ there are exactly $\frac 12 (a-1)$ values $y_j$ satisfying
$T'_a(y_j)=0, T_a(y_j) = T_b(x_i)$.
Hence the number of singular points is $ \frac 12 (a-1)(b-1),$ and they form
two rectangular grids. Since the roots of $T'_b(x)= 0$ are simple,
we see that these points  are crossing points.
\EPf
\pn
\begin{remark}
It follows from their definitions
that $\abs{R} = \frac 12 \pent{b-1}2 (a-1)$,
$\abs{R'} = \frac 12 \pent{b}2 (a-1)$ where $\lfloor x \rfloor$ is the
greatest integer less than or equal to $x$.
\end{remark}
\begin{proposition}
Let $a$ and $b$ are nonnegative coprime integers, a being odd. Let
the Chebyshev curve $\cC$ be defined by the equation  $ T_b(x) -
T_a(y) = 0.$ Then $\cC$ admits the parametrization $ x= T_a(t), \  y=T_b(t).$
The pairs $(t,s)$ giving a crossing point are
$$t=\cos \left ( \Frac ka + \Frac hb
\right ) \pi,  \
s=\cos \left ( \Frac ka - \Frac hb \right ) \pi, $$
where $k,h$ are positive integers such that $\Frac ka + \Frac hb <1.$
\end{proposition}
\Pf
Since $T_a \circ T_b = T_b \circ
T_a=T_{ab}$,  the rational curve $\cC'$ parametrized by $x=T_a(t), \
y=T_b(t)$ is contained in ${\cC}.$
These two curves intersect the line $\{x= x_0\}$ in one point if $\abs{x_0}> 1$, in
$a$ points if $ \abs{x_0}<1$  and in $\frac 12 (a+1)$ points if $x_0= \pm 1.$
Consequently, they are equal.
\pn
The $\frac 12 (a-1)(b-1)$ pairs
$$ t = \cos \left ({k \over a } + {h \over b } \right ) \pi,  \
s= \cos \left ({k \over a } - { h \over b } \right )\pi,  $$
give rise to double  points of $\cC'=\cC$. Because the number of
singular points of $\cC$ is $\frac 12 (a-1)(b-1),$ we see that
there is no other singular point.
\EPf
\begin{remark}
We observe  that the crossing points are obtained for the $(a-1)(b-1)$
elements of
\[
E = \{t_u = \cos \frac{u}{ab} \pi,
\,  0 \leq u \leq ab, \ a \not{\vert} \, u, \ b \not{\vert} \ u \}.\label{E}
\]
For these values, we get $T_b(x(t_u)) = T_a(y(t_u))=(-1)^u$.
Note that $t_u$ and $t_{u'}$ correspond to the same point when
$u \equiv -u' \ \Mod{2b}$ and $u \equiv u' \ \Mod{2a}$.
\end{remark}
\begin{remark}
In general, the curve
$\cC: \ T_b(x) - T_a(y) = 0$ has $\pent d2 +1$ components where
$d=\gcd(a,b)$. See Figure \ref{T510}.
\end{remark}
The following proposition will be useful to consider
Chebyshev curves  as trajectories in a rectangular billiard
(see \cite{JP}).
\begin{proposition}
Let ${\cC}$ be the Chebyshev curve: $T_b(x) - T_a (y)= 0.$
There exists an homeomorphism from the square $I^2=[-1,1]^2$
to the rectangle $ [0,b] \times [0,a]$,
such that the image of $\cC \bigcap I^2 $ is the union of
all the billiard trajectories with slopes $\pm 1$ through
the points with coordinates $ x= b, \  y = a-2k,   0 \le 2k \le a.$
\end{proposition}
\Pf Consider the mapping $F(x,y)=(X,Y)$ with
 $\pi X = b (\pi - \arccos x),  \
\pi Y = a (\pi - \arccos y) .$ By  trigonometry, it is
not hard to check that $F$ has the  announced properties.
\EPf
\begin{figure}[th]
\begin{center}
{\scalebox{.8}{\includegraphics{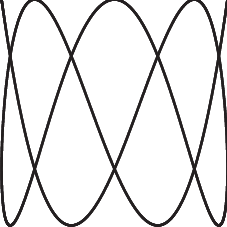}}}\quad\quad
{\scalebox{.7}{\includegraphics{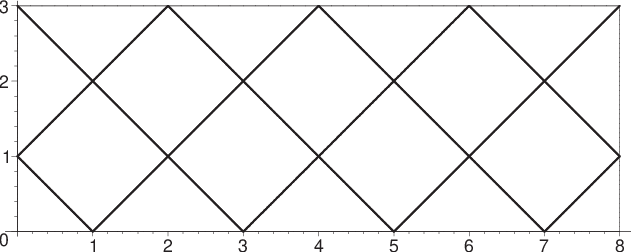}}}
\caption{$T_8(x)=T_3(y)$ and its billiard picture}
\end{center}
\end{figure}
\pn
We shall now present a description of Chebyshev curves using braids.
Let $B_n$ be the group of braids on $n$ strings. For practical
purposes we shall draw these braids horizontally, the strings
being numbered from the bottom to  the top. The standard braid
generators are denoted $\sigma _1, \sigma _2, \ldots,  \sigma
_{n-1}.$ The braid $\sigma_i$ exchanges the strings $i+1$ and $i$,
the string $i+1$ passing over the string $i.$
In this paragraph we shall be interested in plane projections of
braids, called plane braids. We shall also consider the
composition of such plane braids. Let $a=n.$
Let $s_i$ denote the plane braid which is the plane projection of
$\sigma_i$. This plane braid has one crossing point.
\newcommand\seven{s_{\scriptstyle{\rm even}}}
\newcommand\sodd{s_{\scriptstyle{\rm odd}}}

Following Hoste and Zirbel \cite{HZ}, let us define the plane braids $\seven$
and $\sodd $ as
$$ \seven = s_2 s_4 \cdots s_E,  \quad   \sodd =  s_1 s_3
\cdots s_O , $$
where $E$ and $O$ are the largest even and odd integers less than $n=a.$
\begin{figure}[th]
\begin{center}
{\scalebox{.6}{\includegraphics{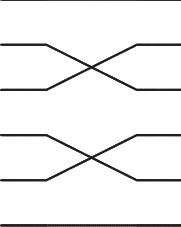}}}\quad
{\scalebox{.6}{\includegraphics{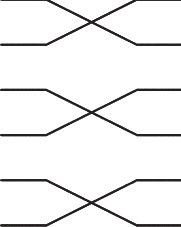}}}\quad \quad      \quad \quad
{\scalebox{.6}{\includegraphics{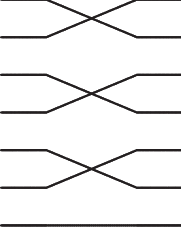}}}\quad
{\scalebox{.6}{\includegraphics{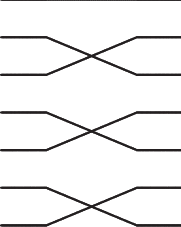}}}
\caption{$\seven$ and $\sodd$ for $n=6, 7$.\label{sevensodd}}
\end{center}
\end{figure}

\begin{proposition}
Let $a,b$ be integers, $a$ being odd. Let ${\cC}$ be the Chebyshev curve
$T_b(x)-T_a  (y)= 0 $. Let $\eps>0$ small enough and consider the
rectangle $ R_{\eps}=\{ |x|<1-\epsilon,  \ |y| \le  1 \}$.
 Then there is a homeomorphism between the pairs
 $(R_{\eps},\cC )$ and $(R_{\eps}, \rho)$ where
 $\rho = (\sodd \, \seven)^{b-1 \over 2 }$
if $b $ is odd and $\rho = (\seven \, \sodd)^{b-2 \over 2} \,
\seven$ if $b$ is even.
\end{proposition}
\Pf Following the proof of Proposition \ref{cc}, the $\frac 12 (a-1)(b-1)$ singular
points of $\cC$ are in $R_\eps$ when $\eps$ is small enough. 
For each $k = 1,\ldots,b-1$, there are $\frac 12 (a-1)$ singular points
$$
(x_k,y_l) = \left (
\cos  k \Frac{\pi}{b},
\cos l \Frac{\pi}{a}
\right ), \ k+l \equiv 0 \ \Mod 2.
$$
It means that over a neighborhood over $x_k$,
the curve ${\cC} $ is isotopic to $\seven$ if $k$ is odd, and isotopic to
$s_{odd} $ if $k$ is even. This proves the result. \EPf
\pn
We can
define the plat closure of a plane horizontal braid with $2m$
strings labelled $ 0, 1, \ldots,  2m-1, $ to be  the plane curve
obtained by connecting the right ends $0$ to $1$, \ldots, $2m-2$
to $2m-1,$ and the left ends in the same order.
\begin{corollary}
Let $a$ be an odd integer, and $b$ an even integer. Let $\rho'$
be the plane braid with $a+1$ strings obtained by adding a free
string numbered $a+1$ over $ \rho = (\seven \, \sodd)^{b-2
\over 2 } \, \seven  .$ Then the Chebyshev curve $T_b(x)-T_a(y)=
0$ is isotopic (in $\SS^2$) to the plat closure of the plane braid $
\rho'.$
\end{corollary}
\Pf Let us illustrate this by looking  at the curve $ T_{10}
(x)-T_5(y)= 0$, which has 3 components. We see on Figure \ref{T510} that it is
the plat closure of $ (\seven\ \sodd)^4 \seven .$ \EPf

\section{Harmonic knots}\label{harmonic}
\begin{figure}[th]
\begin{center}
\begin{tabular}{ccc}
{\scalebox{.2}{\includegraphics{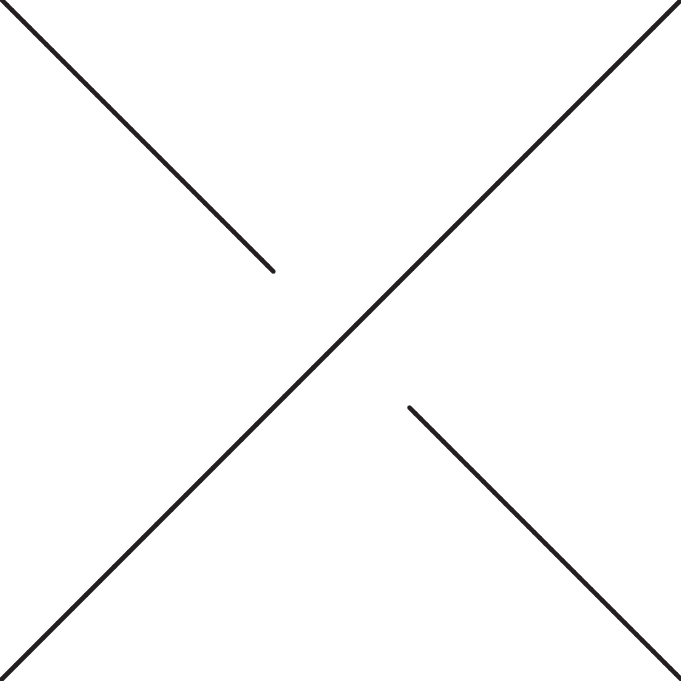}}} &\quad&
{\scalebox{.2}{\includegraphics{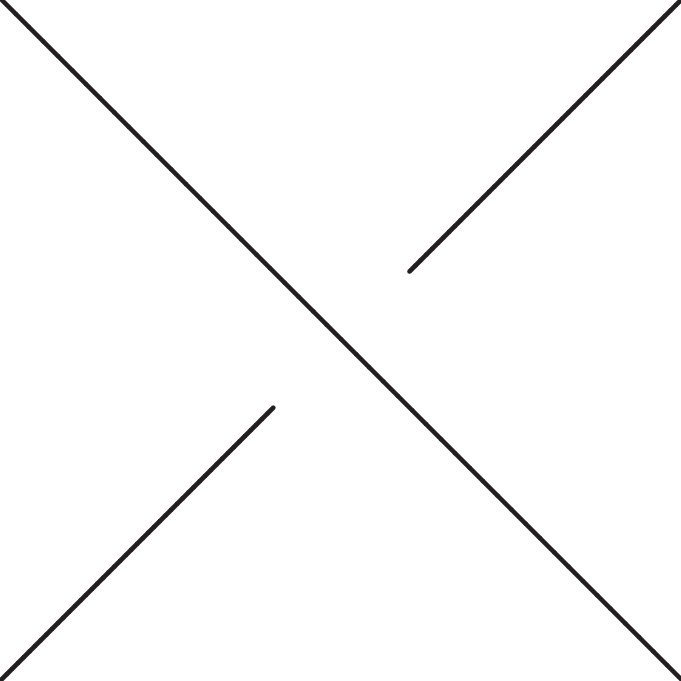}}}\\
$D>0$&&$D<0$
\end{tabular}
\caption{The right twist and the left twist\label{signf}}
\end{center}
\end{figure}

In this paragraph we shall study Chebyshev knots with $ \phi= 0.$
Comstock (1897) found the number of crossing points of the
harmonic  curve parametrized by $x=T_a(t), y=T_b(t), z=T_c(t).$
In particular, he proved that this curve is non singular if and
only if $ a,b,c$ are pairwise coprime integers \cite{Com}. Such
curves will be named harmonic knots  $\H(a,  b,c)$.
\pn
We see that $\H(a,b,1)$ is the unknot because the height function
is monotonic.
We can also obtain the unknot in a less trivial way.
\begin{figure}[th]
\begin{center}
{\scalebox{.4}{\includegraphics{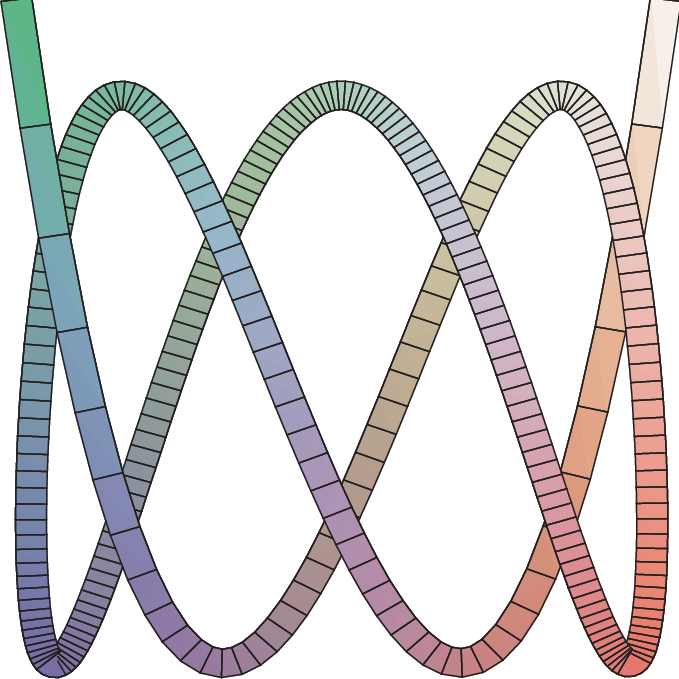}}}
\caption{ The knot $ \H(3,8,11)$ is  trivial. \label{H3811}}
\end{center}
\end{figure}
\begin{proposition}
Let $a,b$ be coprime integers and $c= a+b.$ The harmonic knot
$\H(a,b,c)$ is trivial.
\end{proposition}
\Pf Let $ t \in [-1,1] , \  t=\cos \theta .$ We have
$x=T_a(t)=\cos a \,\theta,  y = \cos b \,\theta, z = \cos  (a+b)\,\theta.$
By trigonometry, we see that the bounded part of our knot is on the surface
$$
S = \{ (x,y,z) \in \RR ^3 ,  \  |x| \le 1, |y| \le 1,z=xy \pm \sqrt{
  (1-x^2) (1-y^2)} \}.
$$
Since $S$  is the union of two sheets that are homeomorphic to  the
square $ [0, 1]^2$ glued along their boundaries, we see that
it is  homeomorphic to a sphere.
Consequently the genus of $\H(a,b,c)$ is zero,
hence it is the  unknot.
\EPf
\pn
Note that the surface $S$ has the symmetries of a regular tetrahedron.
It is contained in the cubic surface
$ \{ x^2+y^2 +z^2 = 1+ 2xyz \} $ which has the same symmetries.
\pn
Let $\cC$ be a plane projection  of a parametrized knot. Consider a crossing
point of $\cC$ obtained for the parameter pair
$(t,s)$. The tangents at this point have opposite slopes (see \cite{KPR}, Lemma 4).
It follows easily that the  nature of this crossing point
depends only of  the sign of the expression
$D=\Bigl(z(t)-z(s)\Bigr)x'(t)y'(t)$. This is not the usual definition
of the sign of oriented crossings, see Figure \ref{signf}.
\begin{lemma}\label{signl}
Let $\H(a,b,c)$ be a harmonic knot.
The nature of the crossing point of parameter
$
t = \cos \left ({k \over a } + {h \over b } \right ) \pi,  \
$
 is given by
$$ {\rm sign}(D) = {\rm sign } \Bigl( -(-1)^{h+k}
\sin( \Frac{ah}b \pi) \, \sin( \Frac{bk}a \pi)
\sin \bigl({ch \over b} \pi \bigr) \,  \sin \bigl({ck \over a} \pi \bigr) \Bigr) .
$$
\end{lemma}
\Pf
The crossing points of the plane projection $ \{ x=T_a(t), \,
y=T_b(t) \} $ are obtained for the parameters
$t =\cos \tau ,  \  s=\cos \sigma $, where
$$
\tau = \left ( \Frac ka + \Frac hb\right ) \pi,  \
\sigma = \left (\Frac ka - \Frac hb \right )\pi.
$$
Using trigonometry we get
$ x'(t)= a \Frac{\sin a\tau}{\sin \tau} , \
y'(t)= b \Frac{\sin b\tau}{\sin \tau}$, so
\[
x'(t) y'(t) =
(-1)^{h+k}\Frac{ab}{\sin^2 \tau}\,
\,\sin( \Frac{ah}b \pi) \, \sin( \Frac{bk}a \pi).\label{x'y'}
\]
We have also
\[
z(t)-z(s) = T_c(t)-T_c(s) =
-2 \,\sin \bigl({ch \over b} \pi \bigr)\, \sin \bigl({ck \over a} \pi \bigr).
\label{ztzs}
\]
and the announced result.
\EPf
\subsection*{Alternate harmonic knots}
The following theorem is the analogue of  a theorem of Lamm
\cite{La} concerning Lissajous knots.
\begin{theorem}[Alternate harmonic knots]
Let $a,b$ be positive coprime integers, and $c=ab-a-b$. The harmonic
diagram  $\H(a,b,c)$ is alternating.
\end{theorem}
\Pf Using Equation (\ref{ztzs}), we get
\begin{eqnarray*}
z(t)-z(s) &= &
-2 \sin \bigl({ch \over b} \pi \bigr) \sin \bigl({ck \over a} \pi \bigr) \\
&= & -2 \sin \bigl((a-1)h\pi - {ah \over b} \pi\bigr)
\sin \bigl((b-1)k\pi - {bk \over a} \pi\bigr)\\
&= & -2 (-1)^{h(a-1) + k(b-1) } \sin \bigl({ah \over b} \pi\bigr) \sin \bigl({bk
\over a} \pi\bigr) .
\end{eqnarray*}
Using Equation (\ref{x'y'}), we get
$ \sign D =  -(-1)^{ah+bk}.$
The crossing points are obtained for the $(a-1)(b-1)$ elements of
$E = \{t_u = \cos \frac{u}{ab} \pi,
\,  0 \leq u \leq ab, \ a \not{\vert} \, u, \
b \not{\vert} \, u \}$.

Note that $t_u < t_{u-1}$ and that at the crossing point corresponding
to $t_u \in E$ one has $\sign D = -(-1)^u$.

The polynomial $x'(t) y'(t)$ has $(a+b-2)$ simple roots
for $t_u = \cos \frac{u}{ab} \pi$, where $a$ or $b$ divides $u =
1, \ldots, ab-1$. For these parameters, the billiard curve
corresponding to the $(x,y)$-plane projection bounces on a wall.

Three cases may occur because at least one of three consecutive $t_u$
belongs to $E$.
\begin{enumerate}
\item $t_{h+1} \in E$ and $t_{h} \in E$. Then
$\sign{x'(t_h)y'(t_h)}= \sign{x'(t_{h+1})y'(t_{h+1})}$
and since the sign of $D$  changes, we conclude that the
sign of $ z(t)-z(s) $ changes between the 2
consecutive parameters $t_{h+1}$ and $t_{h}$.
\item $t_{h+1} \in E$, $t_{h} \not \in E$, $t_{h-1} \not \in E$.
We have $x'(t)y'(t)=0$ at $t_{h}$ and $t_{h-1}$. For $t_{h-1}<t<t_{h-2}$,
we have  $\sign{x'(t)y'(t)}=\sign{x'(t_{h+1})y'(t_{h+1})}$, so
$\sign{x'(t_{h-2})y'(t_{h-2})}=\sign{x'(t_{h+1})y'(t_{h+1})}$.
Hence we see that the sign of $ z(t)-z(s)$ changes between the 2
consecutive parameters $t_{h+1}$ and $t_{h-2}$.
\item $t_{h+1} \in E$, $t_{h} \not \in E$, $t_{h-1} \in E$.
We have $x'(t)y'(t)=0$ at $t_{h}$, so
$\sign{x'(t_{h-1})y'(t_{h-1})}=-\sign{x'(t_{h+1})y'(t_{h+1})}$.
Hence we see that the sign of $ z(t)-z(s)$ changes between the 2
consecutive parameters $t_{h+1}$ and $t_{h-1}$.
\end{enumerate}
In conclusion, the diagram is alternating.
\EPf
\subsection*{Symmetries and  harmonic knots}
A knot $K$ in $\SS^3$ is strongly $(-)$amphicheiral if there is an involution of
$(\SS^3,K)$ which reverses the orientation of both $\SS^3$ and $K$.
A knot $K$ in ${\bf S}^3$ is strongly reversible (or strongly invertible)
if there is an involution of $(\SS^3,K)$ which preserves the orientation of $\SS^3$ and reverses
the orientation of $K$ (see \cite{Kaw}, pp.~127-128).
\begin{proposition}
The harmonic knot $\H(a,b,c)$ is either strongly $(-)$amphicheiral if
$abc$ is odd, or strongly reversible if $abc$ is even.
\end{proposition}
\Pf It is immediate from the parity of Chebyshev polynomials.
\EPf
\begin{corollary}
There are infinitely many amphicheiral harmonic knots.
There are infinitely many strongly reversible harmonic knots.
\end{corollary}
\Pf Since the harmonic knot $\H(a,b,c), \   c=ab-a-b$ is
alternate, its  crossing number is $\frac 12 (a-1)(b-1).$ From this we
conclude that  there is an infinity of  such knots with $a,b,c$ odd, or
with $abc$ even.
\EPf
\pn
If $\sigma$ is any permutation of $\{a,b,c\}$ then the harmonic knot
 $\H(\sigma(a),\sigma(b),\sigma(c))$ is either $\H(a,b,c)$  if
$\sigma$ is an even permutation  or its mirror image  if $\sigma$ is an odd permutation.
\pn
\begin{proposition}\label{abc1}
Let $a, b$ be coprime integers. $\H(a,b,c)$,
$\H(a,b,2ab-c)$ and $\H(a,b,2ab+c)$ are the same knot.
\end{proposition}
\Pf The expression of $\sign D$ (Equations (\ref{x'y'}) and (\ref{ztzs})) for
a given pair of parameters $(t,s)$ corresponding to crossing
points in the $(x,y)$-plane projection is invariant under the
transformation $c \mapsto c+2ab$ and $c \mapsto 2ab-c$.
\EPf
\pn
We can therefore suppose that $a<b$ and $0<c<ab$ to consider all cases.
\begin{proposition}\label{abc2}
Let $a,b,c$ be relatively prime integers. There exists
$c'$ such that $\H(a,b,c')$ is the mirror image of $\H(a,b,c)$.
\end{proposition}
\Pf
Because $a$ and $b$ are relatively prime, one can write
$c = \alpha a + \beta b$ where $\alpha$ and $\beta$ are integers.
Let us consider $c'=-\alpha a + \beta b$. We have $c' \equiv c \ \Mod{2a}$ and
$c' \equiv -c \ \Mod{2b}$. For any crossing point of the Chebyshev
curve $\cC : T_b(x)=T_a(y)$ corresponding to diagrams of both $\H(a,b,c)$ and
$\H(a,b,c')$, we see that $\sign D$ changes to opposite when $c$ is replaced
by $c'$.
\EPf
\begin{corollary}
$\H(a,b,ab+a-b)$ is the mirror image of the alternate knot $\H(a,b,ab-a-b)$.
\end{corollary}
\begin{corollary}\label{phi} Let $a, b$ be relatively prime integers.
There are at most $\phi(a)\phi(b)$ different harmonic knots $\H(a,b,c)$.
\end{corollary}
\Pf
The number of $c$ in $[1,ab]$ that are relatively prime to $a$ and
$b$ is $\phi(a)\phi(b)$ where $\phi$ is the Euler function.
\EPf
\begin{remark}
Because of Propositions \ref{abc1} and \ref{abc2}, for each
$c$ there exists $c'<ab$ such that $\H(a,b,c')$ is the mirror image of
$\H(a,b,c)$. We have at most
$\frac 12\phi(a)\phi(b)$ different knots when we identify the knots and their mirror images.\\
When $a+1<b$, $\H(a,b,1)$ and its mirror image, $\H(a,b,a+b)$ and $\H(a,b,b-a)$  are
trivial.
\pn
Independently but later, G. and J. Freudenburg  (\cite{FF})
proved the following improvement of Proposition \ref{abc2}:
{\em There is a polynomial automorphism $\Phi$
of $\RR^3$ such that $ \Phi (\H(a,b,c)) =\H(a,c',b)$}. They deduce another proof
that the harmonic knots $\H(a,b,a+b)$ are trivial.
\end{remark}
\subsection*{The simplest alternate harmonic knots }
It is remarkable that for $a= 3$ the curves  are drawn in ``Conway
normal form''  for 2-bridge knots \cite{Mu}. Then their Conway
notation  is $ \H(3,n, 2n-3) = C(1,1,\ldots,1)$
when $n$ is not a multiple of $3.$
Turner \cite{Tu} named these knots Fibonacci knots, because their
determinants are Fibonacci numbers. For $n= 4 $ we obtain the
trefoil, for $n= 5$ the figure-eight, for  $n= 7$ the  $ 6_3$
knot, and for $n= 8$ the $\overline{7}_7$ knot.
\begin{figure}[th]
\begin{center}
{\scalebox{.7}{\includegraphics{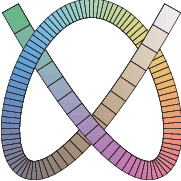}}}   \quad
{\scalebox{.7}{\includegraphics{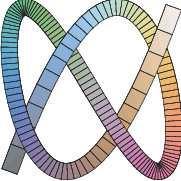}}}   \quad
{\scalebox{.7}{\includegraphics{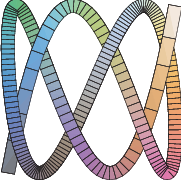}}}
\caption{ The trefoil, the figure-eight knot, and the
$6_3$ knot. \label{AHK}}
\end{center}
\end{figure}
\pn
The Fibonacci knots with an even crossing number are
2-bridge amphicheiral knots.
We have recently proved (\cite{KP5}) that they are not Lissajous.
\pn For $a= 4,$ we also obtain  2-bridge knots. Following the
classical method (\cite{Mu} p. 183-187), we see that their Conway
notation is $\H(4,n, 3n-4) = C (-1,-2,\ldots,-1,-2) =  C(-3,-1,-2,
\ldots,-1,-2),$   ($n$ odd).
\begin{figure}[th]
\begin{center}
{\scalebox{.3}{\includegraphics{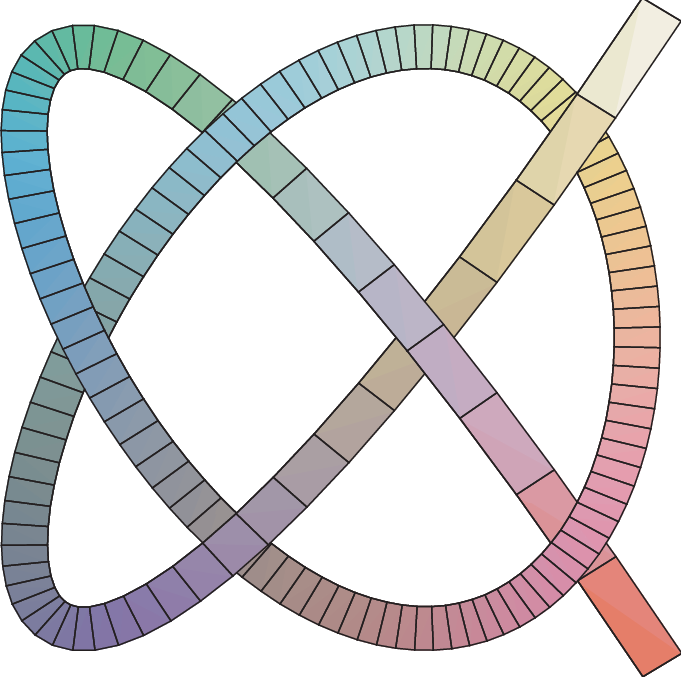}}} \quad
{\scalebox{.3}{\includegraphics{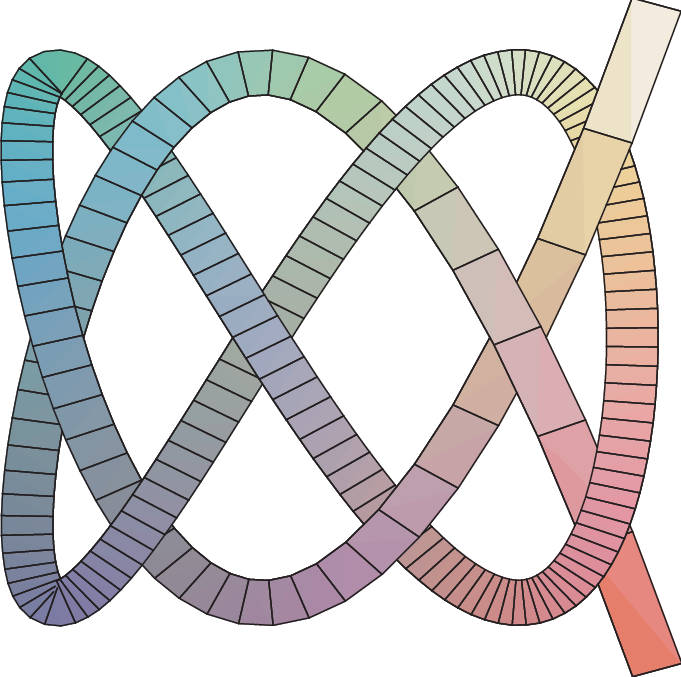}}}
\caption{ The  $6_2$ and $9_{20} $ knots . \label{H4}}
\end{center}
\end{figure}
For $n= 5$ we obtain the  $6_2$ knot,
for $n= 7$ a symmetric picture of the  $9_{20}$ knot (compare
with  Rolfsen's table \cite{Ro}).
\pn
For $a\geq 5$ we may obtain $p$-bridge knots, with $p \geq 3$.
For example the harmonic knot $ \H(5,6,19)$ is the mirror image of $10_{116}$  knot in Rolfsen's table
(amazingly, with exactly the same picture). Its bridge number is known to be $3$.
\begin{figure}[th]
\begin{center}
{\scalebox{.3}{\includegraphics{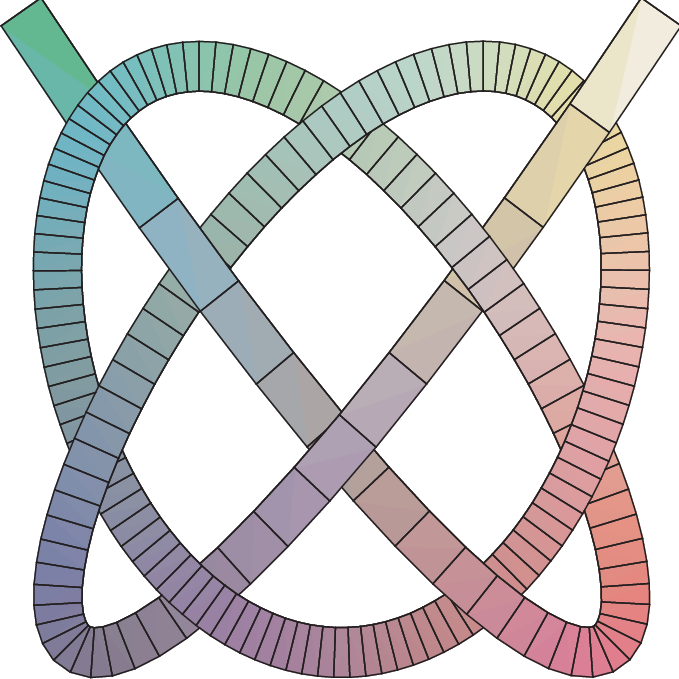}}}
\caption{The mirror image of the $10_{116}$ knot. \label{H5619}}
\end{center}
\end{figure}
\pn
Note that the torus knot $\T(2,2n+1)$ cannot be obtained as an alternate
harmonic knot,
except for the trefoil knot,
as it is proved in \cite{KP1}.
Nevertheless it can be obtained as an harmonic knot.
\subsection*{The torus knots ${\mathbf{\T( 2,2n+1)}}$ }
\begin{theorem}
The knot $\H(3, 3n+2, 3n+1)$ is
the torus knot $\T(2,2n+1)$.
\end{theorem}
\begin{figure}[th]
\begin{center}
{\scalebox{1}{\includegraphics{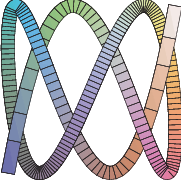}}} \quad
{\scalebox{1}{\includegraphics{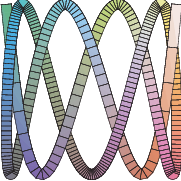}}}
\caption{The torus knots $\overline{5}_1$: $\H(3,7,8)$ and $\overline{7_1}$: $\H(3,10,11)$}
\end{center}
\end{figure}
\Pf
We shall determine the Conway normal form of the harmonic knot
$\H(3,b,c), \    b=3n+1,  c=b+1 .$ The crossing points of the plane projection of
$\H(3,b,c)$ are obtained for pairs of values $(t,s)$
where $ t= \cos \tau , \ s= \cos \sigma,$
and $\tau= \Frac m{3b} \pi , \ \sigma = \Frac {m'}{3b} \pi .$

For $k = 0, \ldots, n-1$, let us consider
\begin{itemize}
\item[] $A_{k}$ be obtained for $m=3 k +1, \  m'= 2b-m.$
\item[] $B_{k}$ be obtained for $ m = 2b+  3 k +2, \ m'= m-2b$.
\item[] $C_{k}$ be obtained for $ m = 2b + 3k + 3, \  m'= 4b-m$.
\end{itemize}
\psfrag{a0}{\small $A_0$}\psfrag{b0}{\small $B_0$}\psfrag{c0}{\small $C_0$}%
\psfrag{a1}{\small $A_1$}%
\psfrag{an}{\small $A_{n-1}$}\psfrag{bn}{\small $B_{n-1}$}\psfrag{cn}{\small $C_{n-1}$}%
\begin{figure}[th]
\begin{center}
{\scalebox{.7}{\includegraphics{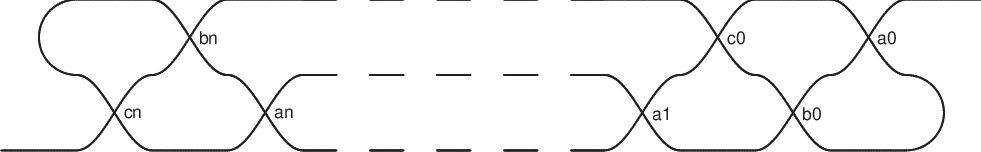}}}
\end{center}
\caption{$\H(3,3n+1,c)$, $n$ even} \label{dh3}
\end{figure}
Then we have
\begin{itemize}
\item[] $x(A_{k}) = \cos ( \Frac {3 k +1}b \pi )$,
$y(A_k) = \frac 12 (-1)^k$.
\item[] $x(B_{k} )= \cos ( \Frac {3 k +2}b \pi )$,
$y(B_k) = \frac 12 (-1)^{k+1}$.
\item[]
$x(C_{k}) = \cos ( \Frac {3 k +3}b \pi )$,
$y(C_k) = \frac 12 (-1)^{k}$.
\end{itemize}
Hence our $3n$ points satisfy
$$x(A_{k-1}) > x(B_{k-1})
>x(C_{k-1}) >  x( A_{k}) > x(B_{k}) > x ( C_{k} ),\ k=1, \ldots , n-1 .$$
Let us determine the nature of the crossing points.
Using the identity $T'_a( \cos \tau ) =  a  \Frac{ \sin a \tau}{\sin \tau },$
we get
$$
\sign{x'(t)y'(t)}= \sign{\sin 3\tau \sin b\tau}.
$$
We get
\begin{itemize}
\item[]
for $A_k$:
$
\begin{array}[t]{rcl}
\sign{x'(t)y'(t)}&=&
\sign{\sin ( \Frac {3 k +1}b \pi ) \sin ( \Frac {3
k +1}3 \pi )}  = (-1)^{k}.
\end{array}$
\item[] for $B_k$:
$
\begin{array}[t]{rcl}
\sign{x'(t)y'(t)}&=&
\sign{\sin \bigl( \Frac {2b +3 k +2}b \pi \bigr) \sin \bigl( \Frac {2b+3
k +2}3 \pi \bigr)}\\
& =&
(-1)^{k+1} \sign{\sin \bigl(\Frac{3k+2}b \pi \bigr) \sin \Frac {\pi}3}=
(-1)^{k+1}.
\end{array}$
\item[] for $C_k$:
$
\begin{array}[t]{rcl}
\sign{x'(t)y'(t)}&=&
\sign{\sin ( \Frac {2b +3 k +3}b \pi ) \sin ( \Frac {2b+3
k +3}3 \pi )}\\
& =&
(-1)^{k+1} \sign{\sin (\Frac{3k+3}b \pi ) \sin \Frac {2\pi}3}=
(-1)^{k+1}.
\end{array}$
\end{itemize}
Now, let us compute the sign of
$$
T_c(t) - T_c(s)=
-2 \sin ( c \Frac {\tau + \sigma}2 ) \sin ( c \Frac {\tau - \sigma}2 )
= -2 \sin \Bigl( \Frac{c}{6b}(m+m')\pi \Bigr)
\sin \Bigl( \Frac{c}{6b}(m-m')\pi \Bigr)
$$
We have, with $c=b+1=3n+2$,
\begin{itemize}
\item[]
for $A_k$:
$
\begin{array}[t]{rcl}
z(t)-z(s) &=&
-2 \sin c \Frac{\pi}3  \sin \bigl( c \Frac{m-b}{3b} \pi \bigr)\\
&=& - 2 (-1)^n \sin \Frac {2 \pi }3 \sin \bigl(c \Frac{n-k}{b} \pi\bigr) \\
&=& 2 (-1)^n \sin \Frac {2 \pi }3 (-1)^{n-k} \sin \bigl(\Frac{n-k}{b} \pi\bigr)
\end{array}
$\\
so $\sign{z(t)-z(s)}= (-1)^k$.
\item[]
for $B_k$:
$
\begin{array}[t]{rcl}
z(t)-z(s) &=&
-2 \sin \bigl( c \Frac{b+3k+2}{3b} \pi\bigr) \sin ( c \Frac{\pi}{3})\\
&=&
-2 \sin \bigl( ( n+k+1+ \Frac {n+ k +1}b) \pi\bigr)\, (-1)^n \sin
\Frac{2\pi}{3} \\
&=& 2 (-1)^k \sin \bigl(\Frac {n+ k +1}b \pi \bigr)\, \sin \Frac{2\pi}{3}.
\end{array}$\\
so $\sign{z(t)-z(s)} = (-1)^k$.
\item[]
for $C_k$:
$
\begin{array}[t]{rcl}
z(t)-z(s) &=&
-2 \sin \bigl( c \Frac {2 \pi }3 ) \sin ( c \Frac {k+1} {b} \pi \bigr)\\
&=&-2 \sin  \bigl( \Frac {4 \pi } 3 \bigr) \sin \bigl(
(b+1) \Frac {k  +1}b \pi \bigr) \\
&=& -2 \sin  \bigl( \Frac {4 \pi } 3 ) (-1)^{k+1} \sin \bigl(\Frac{k+1}b \pi \bigr)
\end{array}$\\
so $\sign{z(t)-z(s)} = (-1)^{k+1}$.
\end{itemize}
\pn
Collecting these results we finally get
$$ \sign{D(A_k)}=1,
\quad
\sign{D(B_k)}= - 1,
\quad
\sign{D(C_k)}= 1.
$$
The Conway  sequence of signs is then
$$
s(C_{n-1}), s(B_{n-1}), s(A_{n-1}), \ldots, s(A_0)
$$
with $s(C_{k}), s(B_{k}), s(A_{k}) = (-1)^{n-1-k}$.
Consequently the Conway normal form of our knot is
$C(1,1,1, -1,-1,-1, \ldots, (-1)^{n-1},(-1)^{n-1},(-1)^{n-1})$.
Its Schubert fraction is $\Frac{2n+1}{2n}\sim - (2n+1)$
and our knot is the torus knot $\T( 2, 2n+1).$
\EPf
\pn
{\bf Remark:} Note that $\H(3,3n+1,3(3n+1)-1)$ is the mirror image of
$\H(3,3n+1,3n+2)$. See Proposition \ref{abc2}.
\pn
{\bf Remark:}
In \cite{KP2}, we obtained the torus knot $\T(2,2n+1)$ as
an alternate polynomial knot where $x(t)=T_3(t), \, y(t)=P(t), \, z(t)=Q(t)$
are polynomials and $\deg P = 3n+1$, $\deg Q = 3n+2$, that is to say the
same degrees.
\pn
Because of their  definitions, the symmetries of the harmonic
knots are easy to find. They are either strongly negative
amphicheiral if $a,b,c$ are odd, or strongly reversible. So that
not every knot is an harmonic knot.
We can also remark that harmonic knots are billiard knots in a convex
(compact) billiard (in fact a truncated cube)\cite{JP}.

On the other hand, it is not difficult to see that if we change the
nature of one crossing point in the diagram of the $10_{116}$ knot,
we can obtain the $8_{17}$ knot. The knot $8_{17}$ is famous because it
is the first non reversible knot.

In the next paragraph, we shall see that it is possible to choose
the nature of the crossing points with a (shifted) Chebyshev
polynomial as height function.
\section{Every knot is a Chebyshev knot}\label{tk}
Let us denote  $B_n$ the group of $n$-braids and $S_n$ the symmetric
group. The group of pure braids $P_n$ is the kernel of the
morphism
 $\pi \, : \  B_n \rightarrow S_n.  \, $
If $\alpha $ is a braid, we shall denote $\rho (\alpha) $ its
plane projection. In the next theorem, which is analogous to a theorem of Lamm
 for Lissajous curves (see \cite{BDHZ,La2}), we show that
\begin{theorem}
Every knot has a projection which is a Chebyshev plane curve.
\end{theorem}
This is a consequence of the following proposition.
\begin{proposition}
Let $K$ be a knot, br$(K)$ its bridge number. Let $m \ge {\rm br}(K)$ be an integer.
Then $K$ has a projection which is a Chebyshev curve
$ x= T_a(t), \  y= T_b(t), $
where $a= 2m-1,$ and $b \equiv 2 \ \Mod{2a}$.
\end{proposition}
\Pf Let $K$ be a knot.
Let $D$ be a regular diagram of $K$ such that the abscissa has only two extremal values
reached at $m$ maxima and $m$ minima.
It means that $K$ is the plat closure of
a horizontal braid $t$ with $ 2m$ strings.  We can suppose the last
string unbraided.
Furthermore, reordering if necessary
the ordinates of the $2m$ extrema, we can suppose that $ \pi (t)
=(2, 3) \cdots (2m-2, 2m-1)=\pi (\sigma _2 \  \sigma_4 \cdots
\sigma_{2m-2}) .$ Let us denote $ \seven=\rho(\sigma_2
\cdots\sigma_{2m-2}), $ and $ \sodd=\rho(\sigma_1 \cdots
\sigma_{2m-3}).$

As $ \pi (t)= \pi (\sigma_2 \cdots \sigma_{2m-2})$ we see that
there  exists $ l \in \ker (\pi)=P_{2m-1},$ the group of pure
braids, such that $t = l \cdot \,  \sigma_2 \cdots \sigma_{2m-2} . $
As the braids  $ A_{i \, j}= x^{-1} \sigma _i ^2 x$, where
$x=\sigma_{i+1} \cdots \sigma_j,$  generate $P_{2m-1},$ the braid $l$
is a composition of such elementary braids. It is not difficult to
see that there is a braid $\alpha _{i \, j} $ equivalent to $
A_{i \, j} $  with plane projection
 $ \rho (\alpha _{i \, j}) = (\seven \sodd)^{a} .$
Consequently, the braid $t$ is equivalent to a braid
projecting on $(\seven \, \sodd)^{ka} \seven .$
 Using the braid description of Chebyshev curves (corollary 1),
we conclude that our knot $K$ is equivalent to a knot projecting
upon the Chebyshev curve
 $ x=T_a(t), \ y = T_b(t),    \    a \ {\rm odd} , \ b \equiv 2 \ \Mod{2 a}.$
 \EPf
\begin{figure}[th]
\begin{center}
{\scalebox{.4}{\rotatebox{90}{\includegraphics{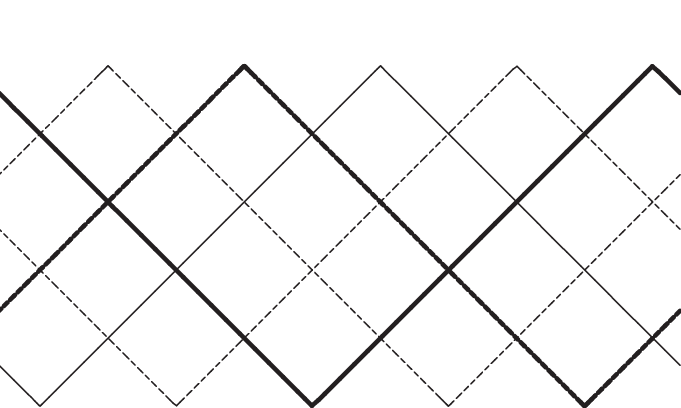}}}}
\caption{The ``plane braid'' $(\seven\, \sodd)^5 $ is
the projection of the braid $\alpha_{2,5}\in B_5.$}
\end{center}
\end{figure}
\pn
We shall prove our principal result with a density argument based on
Kronecker's theorem (\cite{HW}, Theorem 443, p. 382.)
Let us recall this theorem.
\begin{theorem}[Kronecker]
If $ \theta_1, \theta_2, \ldots,  \theta_k,\, 1 $ are linearly
independent over $ {\bf Q},$ then the set of points
$((n\,\theta_1),(n\,\theta_2), \ldots,  (n\,\theta_k))$ is dense
in the unit cube. {\em Here $(x)$ denotes the fractional part of
$x.$}
\end{theorem}
We shall need the following lemma.
\begin{lemma}
Let $c_1,\ldots, c_k$ be real numbers such that $
-1<c_1<c_2<\cdots<c_k<1.$  There exists a positive number $e <1-c_k$
such that the numbers $\arccos (c_1+e), \arccos (c_2 +e),
\ldots,  \arccos (c_k +e),  \, 1 $  are linearly independent over
${\bf Q} .$
\end{lemma}
\Pf First, we shall prove, by induction on $k$, that the functions
$\arccos (c_i +x), \, 1$ are linearly independent over $ {\bf R}.$
Let $ \lambda_0 + \sum _{i= 1} ^k \lambda_i \arccos (c_i +x) = 0$
be a  linear relation  between these functions.
We get by derivating
$$ \sum_{i= 1} ^k \lambda_i { 1 \over \sqrt {1- (c_i +x)^2 }} \, \,=0 .$$
Then, when $ x \rightarrow 1-c_k, $   we get $ \lambda_k=
0,$ and the result follows by induction.
\pn
Suppose now that for
each $e<1-c_k$ there exists a relation
$$\lambda_0 + \sum_{i= 1}^k \lambda_i \arccos (c_i + e) = 0,  \quad
{} \
  \lambda_i \in{\bf Q}, \quad {} \
\sum_{i=1}^k |\lambda _i | = 1 .$$ By cardinality, there are infinitely many
$e$ in $(0, 1-c_k)$ with the same collection of $\lambda_i.$ This
means that the analytic  function
$\lambda (x) = \lambda_0 + \Sum_{i= 1}^k \lambda_i \arccos (c_i  +x)$
has an infinity of zeroes in the interval
 $(0, 1-c_k)$, which is absurd.
\EPf
\begin{theorem}
Every knot is a Chebyshev knot.
\end{theorem}
\Pf Let $K$ be a knot projecting on the Chebyshev curve
$\{x= T_a(t), \ y= T_b(t)\}.$ The crossing points of the projection are
obtained for the distinct pairs of values
$$t=\cos \left ( \Frac ka + \Frac hb
\right ) \pi,  \
s=\cos \left ( \Frac ka - \Frac hb \right ) \pi, \
\Frac ka + \Frac hb <1.$$
Let us denote these values by $(t_i, s_i), \ i=
1\ldots n= \frac 12 (a-1)(b-1).$ By our lemma,
let $e<1-\cos \Frac{\pi}{ab}$ be a positive  number
such that the $2n+1$ numbers
$ \  1, \ \tau_i = \arccos(t_i+e), \, \sigma_i= \arccos (s_i +e),  \ i= 1 ,
\ldots n,  \,$  are linearly independent over ${\bf Q}.$ Let
us define the function $ Z(t)=T_c(t+e)$ (depending on the integer
$c$). We have
$$ Z(t_i)-Z(s_i)=\cos c\,\tau_i - \cos c \,\sigma_i
.$$
Since the numbers $ \, 1,\tau_i, \sigma_i$ are linearly independent over ${\bf Q},$ the
numbers $ c\, \tau_i  \ \Mod{2 \pi}$ and $c\, \sigma_i \ \Mod{2\pi}$ are
dense in $ [0, \, 2\pi ]^{2n}$  by Kronecker's theorem. So that we
can choose arbitrarily the signs of $Z(t_i) - Z(s_i), $ that is,
the over/under nature of the crossing points.
\EPf
\subsubsection*{Example 1: the knot ${\mathbf 6_1}$}
Let us consider the curve parametrized by
$$
x = T_3(t), \, y=T_8(t), \, z = T_{10}(t+\frac{1}{100}).
$$
\begin{figure}[th]
\begin{center}
{\scalebox{.4}{\includegraphics{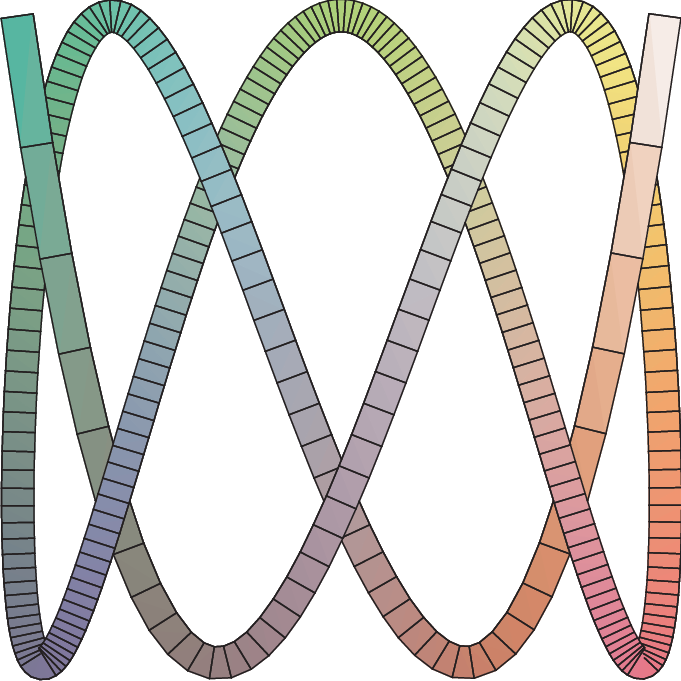}}}
\caption{The  knot $6_1$ is a Chebyshev knot \label{K61}}
\end{center}
\end{figure}
Computing $\sign D$ for the $\frac 12 (3-1)(8-1) = 7$ crossing points, we
find the Conway normal form: $[-1,-1,-1,-1,1,1,1]$. Its Schubert fraction
(see \cite{Mu}) is then
$$
-1 + \Frac{1}{-1 + \Frac{1}{-1+\Frac{1}{-1+\Frac{1}{1+\Frac{1}{1+\Frac 11}}}}} =
\Frac{9}{-5} \sim \Frac{9}{4}.
$$
This knot is the knot $6_1$.
\subsubsection*{Example 2: the knot ${\mathbf 8_{17}}$}
The famous $8_{17}$ knot is non reversible and strongly ($-$)amphicheiral
(see \cite{Kaw} p. 128).
It is  a 3-bridge knot. The Chebyshev curve $T_{6}(x)=T_5(y)$ is
 one of its diagrams.
\begin{figure}[th]
\begin{center}
{\scalebox{.4}{\includegraphics{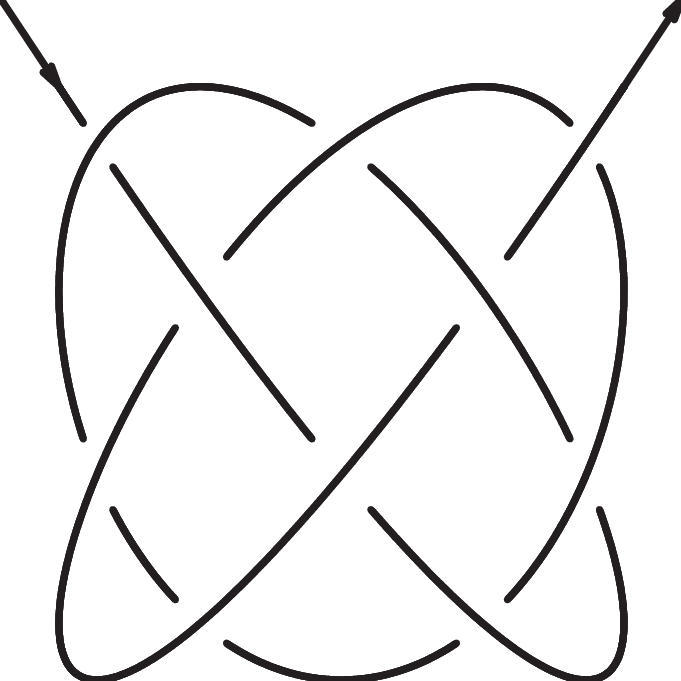}}} \quad
{\scalebox{.4}{\includegraphics{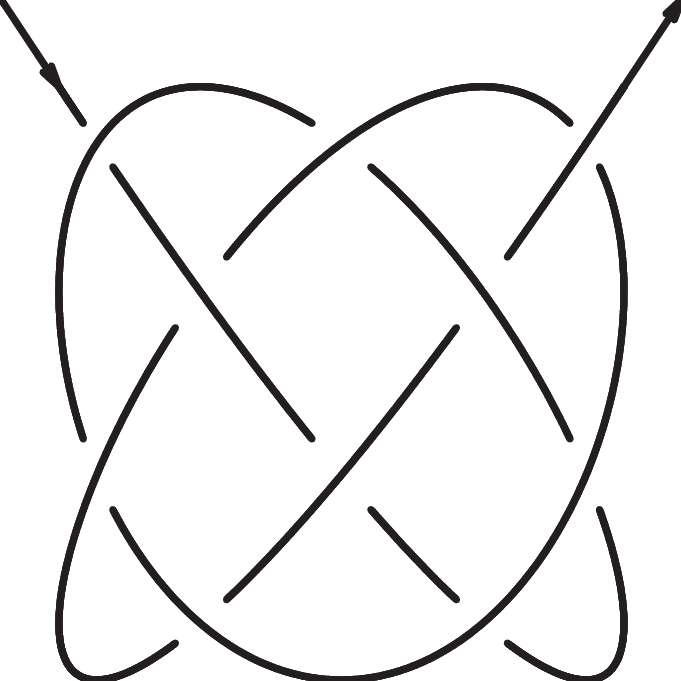}}} \quad
\caption{The  knot ${8}_{17}$ and its reverse as Chebyshev knots \label{K817}}
\end{center}
\end{figure}
It can be parametrized by
$x=T_5(t), y = T_6(t), z= T_{33}(t+148 \cdot 10^{-4})$.
Its reverse can be parametrized by
$x=T_5(t), y = T_6(t), z= T_{33}(t-148 \cdot 10^{-4})$.
We see that it is the reversed image of $8_{17}$
by a half-turn about the $y$ axis.
\section*{Conclusion}
Let us give a list of the first 2-bridge harmonic knots with
their Conway-Rolfsen numbering.
Because of their simplicity, we also give their Chebyshev diagrams.
A bar over a knot name indicates mirror image.
\begin{center}
\begin{tabular}{cccc}
{\scalebox{.7}{\includegraphics{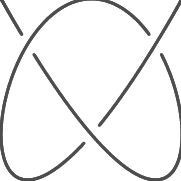}}}&
{\scalebox{.7}{\includegraphics{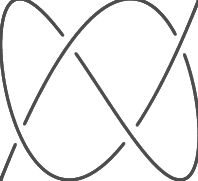}}}&
{\scalebox{.7}{\includegraphics{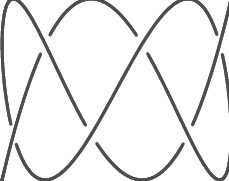}}}&
{\scalebox{.7}{\includegraphics{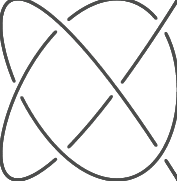}}}\\
$\overline{3}_1$& $4_1$ &$\overline{5}_1$ & $\overline{5}_2$\\
$\H(3,4,5)$&$\H(3,5,7)$&$\H(3,7,8)$&$\H(4,5,7)$
\end{tabular}
\pn
\begin{tabular}{cccc}
{\scalebox{.7}{\includegraphics{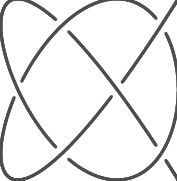}}}&
{\scalebox{.7}{\includegraphics{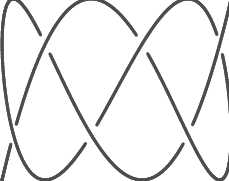}}}&
{\scalebox{.7}{\includegraphics{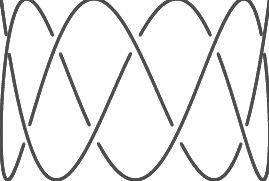}}}&
{\scalebox{.7}{\includegraphics{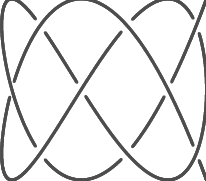}}}\\
$6_2$&$6_3$&$\overline{7}_1$& $7_5$\\
$\H(4,5,11)$&$\H(3,7,11)$&$\H(3,10,11)$&$\H(4,7,9)$
\end{tabular}
\pn
\begin{tabular}{cccc}
{\scalebox{.7}{\includegraphics{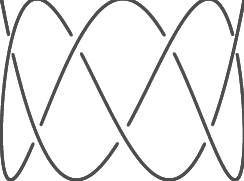}}}&
{\scalebox{.7}{\includegraphics{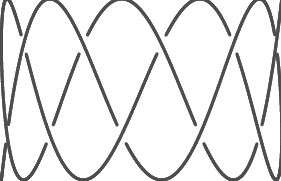}}}&
{\scalebox{.7}{\includegraphics{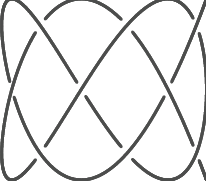}}}&
{\scalebox{.7}{\includegraphics{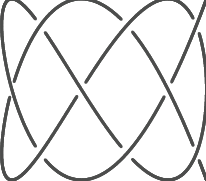}}}\\
$\overline{7}_7$&$8_3$&$8_7$&$9_{20}$\\
$\H(3,8,13)$&$\H(3,11,13)$&$\H(4,7,13)$&$\H(4,7,17)$\\
\end{tabular}
\pn
\begin{tabular}{ccc}
{\scalebox{.7}{\includegraphics{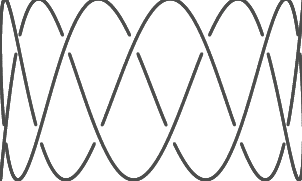}}}&
{\scalebox{.7}{\includegraphics{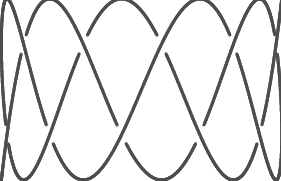}}}&
{\scalebox{.7}{\includegraphics{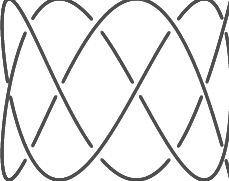}}}\\
$\overline{9}_1$ &$\overline{9}_{17}$&$\overline{9}_{18}$\\
$\H(3,13,14)$&$\H(3,11,16)$&$\H(4,9,11)$
\end{tabular}
\pn
\begin{tabular}{ccc}
{\scalebox{.7}{\includegraphics{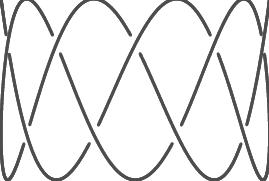}}}&
{\scalebox{.7}{\includegraphics{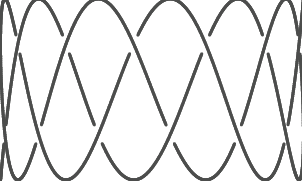}}}&
{\scalebox{.7}{\includegraphics{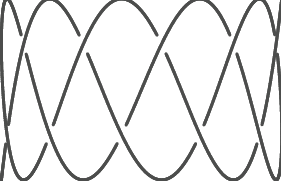}}}\\
$\overline{9}_{31}$ &$10_{37}$&$10_{45}$\\
$\H(3,10,17)$&$\H(3,13,17)$&$\H(3,11,19)$
\end{tabular}
\end{center}
In \cite{KPR}, we give a complete list of Chebyshev parametrizations of
the 2-bridge knots of 10 crossings or less.
\pn
Now, let us give the list of harmonic knots $\H(5,6,c)$.
We get $4=\frac 12 \phi(5) \phi(6)$ different types
up to mirror symmetry.
\begin{center}
\begin{tabular}{cccc}
{\scalebox{.8}{\includegraphics{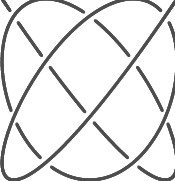}}}&
{\scalebox{.8}{\includegraphics{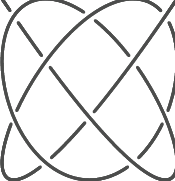}}}&
{\scalebox{.8}{\includegraphics{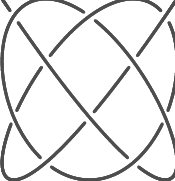}}}&
{\scalebox{.8}{\includegraphics{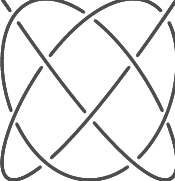}}}\\
$0$ & ${5}_2$& $10_{159}$ & $\overline{10}_{116}$ \\
$\H(5,6,1)$&$\H(5,6,7)$&$\H(5,6,13)$&$\H(5,6,19)$
\end{tabular}
\end{center}
The $10_{159}$ knot is the harmonic knot $ \H(5,6,13)$. Its bridge number is equal to $3$. Our
Chebyshev parametrization provides an easy proof that it is strongly reversible (compare \cite{Kaw},
Appendix F, p. 254).

In conclusion, we have found a great number of distinct harmonic knots. Furthermore,
their diagrams have a small number of crossing points. We hope that our
Chebyshev models will be useful for the study of knots.

In \cite{KP4}, we classify the harmonic knots $\H(3,b,c)$ and
$\H(4,b,c)$. Even for $a= 5$,
the classification of harmonic knots seems to be a difficult and
interesting problem.

\end{document}